\documentclass[11pt]{amsart}
\usepackage{amssymb}
\numberwithin{equation}{section}

\newtheorem{theorem}{Theorem}[section]

\newtheorem{proposition}[theorem]{Proposition}
\newtheorem{lemma}[theorem]{Lemma}
\theoremstyle{definition}
\newtheorem{definition}[theorem]{Definition}
\theoremstyle{remark}
\newtheorem{remark}[theorem]{Remark}
\newtheorem{example}[theorem]{Example}
\newtheorem{notation}[theorem]{Notation}
\newcommand{\U}{\mathcal{U}}

\newcommand{\R}{\mathbb{R}}

\newcommand\lie[1]{\mathfrak{#1}}

\newcommand{\fh}{\lie{h}}
\newcommand{\fg}{\lie{g}}

\newcommand{\fk}{\lie{k}}

\newcommand{\fn}{\lie{n}}

\newcommand{\f}{\frac}

\newcommand{\lam}{\lambda}

\renewcommand{\l}{\langle}
\renewcommand{\r}{\rangle}

\newcommand{\inv}{^{-1}}
\renewcommand{\ker}{ \operatorname {ker}}
                                     
\begin{document}

\title{Relative equilibria at singular points of the moment map }

\author{E. Lerman}
\address{Department of
Mathematics, University of Illinois, Urbana, IL 61801}
\email{lerman@math.uiuc.edu}
\date{\today}

\begin{abstract}
We prove a criterion for stability of relative equilibria in symmetric
Hamiltonian systems at singular points of the momentum map.  This
generalizes a theorem of G.W. Patrick.  The method of the proof is
also useful in studying the bifurcation of relative equilibria.
\end{abstract}

\maketitle

\tableofcontents

\section{Introduction}

Let $(M, \omega)$ be a symplectic manifold, $G$ a Lie group acting
properly on $M$ and preserving the symplectic form.  Assume further
that there is an equivariant moment map $\Phi:M \to \fg^*$
corresponding to the action.  Let $h$ be a $G$ invariant smooth
function on $M$.  Denote its Hamiltonian vector field by $X_h$ and the
flow of the vector field by $e^{tX_h}$.  Since the function $h$ and
the form $\omega$ are invariant, the flow $e^{tX_h}$ is $G$
equivariant.  Hence it descends to a flow on the quotient space $M/G$
(throughout the paper we will assume for simplicity that all flows
exist for all times).

Recall that a point $m$ in $M$ is a {\bf relative equilibrium} of the
Hamiltonian $h$ (really of the flow $e^{tX_h}$) if the orbit $G\cdot
m$, thought of as a point in $M/G$, is fixed by the induced flow on
$M/G$.\footnote{``$\cdot$'' denotes the action of $G$ on $M$}
Equivalently, $m$ is a relative equilibrium iff
\begin{itemize}
\item the vector $X_h(m)$ is tangent to the orbit $G\cdot m$ iff 
\item there is $\xi \in \fg$, the Lie algebra of $G$, such that the
corresponding vector field $\xi_M$ on $M$ satisfies $\xi_M (m) = X_h (m)$ iff
\item $m$ is a critical point of $h -\l \Phi, \xi\r$ for some $\xi \in
\fg$.
\end{itemize}
\begin{remark}\label{rk_isotropy}
It is standard that $\xi_M (m) = X_h (m)$ implies that $e^{t\xi} \cdot
m = e^{tX_h} (m)$ for all $t$.  Since the moment map $\Phi$ is
constant along the flow of $X_h$ and is equivariant, it follows that
$\Phi (m) = \Phi(e^{tX_h} (m)) = \Phi (e^{t\xi} \cdot m ) = e^{t\xi}
\cdot \Phi (m)$.  Hence $\xi$ is in the isotropy Lie algebra of $\Phi
(m)$.
\end{remark}

Since we are dealing with a Hamiltonian system, it does not make sense
to talk about asymptotic stability.  The best one can hope for is
something along the lines of:
\begin{definition}
Let $(M,\omega, G, \Phi, h)$ be as above --- a symplectic manifold
$(M, \omega)$ with a symmetry group $G$ acting properly on $M$ in a
Hamiltonian fashion with a moment map $\Phi$ and with a $G$-invariant
Hamiltonian $h$. Let $H$ be a subgroup of $G$.  A relative equilibrium
$m$ of $h$ is {\bf $H$-stable} if for any $H$ invariant neighborhood
$V$ of the orbit $H\cdot m$ there is an $H$ invariant neighborhood
$U\subset V$ such that $e^{tX_h} (U) \subset V$ for all times $t$.
\end{definition}
\noindent
In his thesis \cite{Patrick}, G. Patrick proved the following result
(see also \cite{Patrick2}):\\[-3pt]

\noindent
{\bf Theorem} [Patrick] {\em Given a Hamiltonian system $(M,\omega, G,
\Phi, h)$ as above, assume additionally that
\begin{enumerate}
\item A relative equilibrium $m$ is a regular point of the moment map
$\Phi$,
\item \label{ass2} the isotropy group $H$ of $\mu = \Phi (m)$ acts
properly on $M$,
\item there exists an $Ad (H)$ invariant inner product on the Lie
algebra $\fg$ of G,
\item \label{ass4} for an element $\xi \in \fg$ with $\xi_M(m) = X_h (m)$ the
Hessian of $h - \l \Phi, \xi\r$ is definite when restricted to a
complement of the tangent space $T_m (H\cdot m)$ to the orbit $H\cdot
m$ in in the tangent space to the level set $T_m (\Phi \inv (m))$.
\end{enumerate}
Then the relative equilibrium $m$ is $H$ stable.}\\[-12pt]

\begin{remark}[Hessians]\label{remark_hessian}
Let $N$ be a manifold, $f\in C^\infty (N)$ a function and $x\in N$ a
point.  The differential $df(x)$ of $f$ at $x$ is a well-defined
functional on the tangent space $T_xX$ which behaves well under
pull-backs: if $\psi : Z \to N$ is a smooth map  of manifolds with
$\psi (z) = x$, then $d(\psi^*f) (z) = \psi^* (df(x))$.  

 On the other hand, unless $df(x) = 0$, the Hessian $d^2 f(x)$ is not
a well-defined quadratic form on $T_xN$.  Of course, for every choice
of coordinates on $N$ near $x$ we get a symmetric matrix of partial
derivatives of $f$.  But this matrix behaves badly under a change of
coordinates and, more generally, under pull-backs by smooth maps.

If $df(x) = 0$ then the Hessian $d^2 f(x)$ is well-defined\footnote{
Recall that $df(x)$ is the class of $f - f(x)$ in $\lie{M}/\lie{M}^2$
where $\lie{M}$ is the ideal of functions which vanish at $x$.  If
$df(x) = 0$, then $f - f(x)\in \lie{M}^2$ and defines $d^2 f(x)$ to be
the class of $f - f(x)$ in $\lie{M}^2/\lie{M}^3$.}  
and behaves well under pull-backs and restrictions: if $\psi : Z \to
N$ is a smooth map of manifolds with $\psi (z) = x$, then
$d^2(\psi^*f) (z) = \psi^* (d^2f(x))$.  In particular, if $Z$ is a
submanifold of $N$, then
\begin{equation}
d^2f(x)|_{T_x Z} = d^2 (f|_Z) (x).
\end{equation}
\end{remark}

\begin{remark}
Note that that since $m$ is assumed to be a regular point of the
moment map, the tangent space to the level set of $\Phi$ at $m$, $T_m
\Phi\inv (m)$, is $\ker (d\Phi _m)$.  Further, a complement to $T_m
(H\cdot m)$ in $T_m \Phi\inv (m)$ is a symplectic subspace isomorphic
to the tangent space of the reduced space $\Phi\inv (\mu)/H$ at the
point $H\cdot m$.  Hence assumption (\ref{ass4}) in the theorem above
amounts to saying that the Hessian of the reduced Hamiltonian is
positive definite at the point $H\cdot m$.  In particular the
equilibrium point $H\cdot m$ is stable in the reduced system.
\end{remark}

A simple example below shows that assumption of stability in the reduced
system is not enough to guarantee that the corresponding relative
equilibrium is relatively stable.

In this paper we extend Patrick's result to the case where the
relative equilibrium $m$ in question is not necessarily a regular
point of the moment map.  Patrick's theorem follows as a special case.
The difficulty with extending Patrick's result lies in the fact that
in general there is no reason for the reduced space $\Phi\inv (\mu)/H$
to be an orbifold in a neighborhood of the point $H\cdot m \in\Phi\inv
(\mu)/H$ (We use the notation of Patrick's Theorem). So while the
Hamiltonian dynamics on the reduced space and the notion of stability
still make sense (\cite{SjL}, \cite{LB}, \cite{ACG}), the notions of
the tangent space and of the Hessian of the reduced Hamiltonian do
not.  The tangent space at $m$ to $\Phi\inv (\mu)$ does not make sense
either (except in the sense of Zariski).

On the other hand the kernel of the differential $d\Phi _m$ is still a
vector space and it still contains $T_m (H\cdot m)$ as a subspace.  So
it makes sense to replace the assumption (\ref{ass4}) of Patrick's
theorem ({\em `` for an element $\xi \in \fg$ with $\xi_M(m) = X_h (m)$ the
Hessian of $h - \l \Phi, \xi\r$ is definite when restricted to a
complement of the tangent space $T_m (H\cdot m)$ to the orbit $H\cdot
m$ in in the tangent space to the level set $T_m (\Phi \inv (m))$''})
by \\[-3pt]

\noindent
{\em for a well-chosen element $\xi \in \fh$ with $\xi_M(m) = X_h (m)$
the restriction of the Hessian of $h - \l \Phi, \xi\r$ to $\ker (d\Phi
_m)$ is semi-definite and the kernel of the restriction is {\bf precisely}
$T_m (H\cdot m)$, the tangent space to the orbit $H\cdot m$.}\\[-3pt]

Since for any $\xi \in \fh$ the restrictions $\l \Phi , \xi\r|_{H\cdot
m}$ is constant and since $h|_{H\cdot m}$ is constant as well, we have
that the restriction $d^2 (h - \l \Phi, \xi\r) (m) |_{T_m (H\cdot m)}$
of the Hessian of $h-\l \Phi, \xi\r)$ to the tangent space to the
$H$-orbit through the relative equilibrium is always zero.  Thus $d^2
(h - \l \Phi, \xi\r) (m) |_{\ker d\Phi (m) }$ descends to a
well-defined quadratic form on $V:= \ker (d\Phi _m)/T_m (H\cdot m)$.

\begin{remark}[symplectic slice] \label{sympl_slice}
Note that since $\ker (d\Phi _m)$ is the symplectic perpendicular $T_m (G\cdot
m)^\omega$ and since $T_m (H\cdot m) = T_m (G\cdot m)^\omega \cap T_m
(G\cdot m)$, it follows that $V:= \ker (d\Phi _m)/T_m (H\cdot m)$ is
also a symplectic vector space.  It is  called the {\bf symplectic
slice} to the orbit $G\cdot m$.  If $m$ is a regular point
of the moment map, then, up to an action of finite group,  the
symplectic slice is symplectomorphic to a neighborhood of the point
$H\cdot m$ the corresponding reduced space.
\end{remark}

Thus the main result of this paper reads:
\begin{theorem}\label{main}
Let $(M, \omega)$ be a symplectic manifold, $G$ a Lie group acting
on $M$ and preserving the symplectic form.  Assume further
that there is an equivariant moment map $\Phi:M \to \fg^*$
corresponding to the action.  Let $h$ be a $G$ invariant smooth
function on $M$.  Suppose $m$ is a relative equilibrium of $h$.  Suppose
\begin{enumerate}
\item the isotropy group $H$ of $\mu = \Phi (m)$ acts properly on $M$,
\item there exists an $Ad (H)$ invariant inner product on the Lie
algebra $\fg$ of $G$,
\item the restriction of the Hessian of $h - \l \Phi, \xi\r$ to $\ker
(d\Phi _m)$ descends to a definite form on the symplectic slice at
$m$, where $\xi \in \fh$ is orthogonal (with respect to the $Ad
(H)$-invariant inner product) to the isotropy Lie algebra $\fg_m$ of
$m$ and satisfies $\xi_M(m) = X_h (m)$.
\end{enumerate}
Then $m$ is $H$ stable.
\end{theorem}
\begin{remark}
Since Theorem~\ref{main} is local, one can relax the assumption on the
action of $H$.  Namely, it is enough to assume that $H$ acts properly
on a neighborhood of the orbit $H\cdot m$.
\end{remark}

We end the section with an example illustrating that a relative
equilibrium may be relatively stable within a  level set, and  
unstable in the whole phase space.
\begin{example}\label{the_example}
Consider the standard action of $SO(2)$ on $\R^2$.  It lifts to a
Hamiltonian action on $T^* \R^2 =\{(q_1, q_2, p_1, p _2)\}$ with a
moment map $f(q, p) = q_1 p_2 - q_2 p_1$.  Now consider a Hamiltonian
system on the product $M = T^* \R^2 \times T^* S^1 = \{(q_1, q_2, p_1,
p _2, \theta, p_\theta )\}$ (standard symplectic form) with the
Hamiltonian
$$
h = (q_1 p_2 - q_2 p_1) + p_\theta (p_1^2 + p_2^2 - q_1 ^2 - q_2 ^2).
$$
The Hamiltonian $h$ Poisson commutes with $p_\theta$, which is a
moment map for the action of $S^1$ on the second factor.

The reduced spaces $M_\mu$ are all symplectomorphic to $T^* \R^2$ and
the reduced Hamiltonians $h_\mu$ are given by
$$
h_\mu = (q_1 p_2 - q_2 p_1) + \mu(p_1^2 + p_2^2 - q_1 ^2 - q_2 ^2).
$$
Note that the origin $(0,0,0,0) \in T^* \R^2$ is stable for $h_0 = q_1
p_2 - q_2 p_1$ and that the Hessian of $h_0$ is not positive definite
at the origin.  For $\mu \not = 0$, the reduced Hamiltonian is of the
form $h_\mu = f + \mu g$ with $\{f, g\} = 0$.  It is easy to see that
the origin is unstable for the flow of $g$, hence is unstable for
$h_\mu$ for all $\mu \not =0$.  Therefore the relative equilibria of
the form $(q_1, q_2, p_1, p_2, \theta, 0)$ are all unstable in $M$,
even though the corresponding fixed point is stable in the reduced
space $M_0$.
\end{example}

\subsection*{Acknowledgments}
My original motivation for this work was the study of a two-body
problem on $S^2$, a problem proposed by Larry Bates.  Debra Lewis's
critical comments at the AMS meeting in Corvallis have been very
helpful.  I thank Tudor Ratiu and Juan-Pablo Ortega for making
available to me the early versions of their paper {\em Persistance et
diff\'erentiabilit\'e de l'ensemble des \'el\'ements relatifs dans les
syst\`emes hamiltoniens symmetriques} and for making me aware of the
results of Montaldi on the stability of relative equilibria
\cite{Montaldi}.  I thank James Montaldi for promptly responding to my
query.

\section{Reduction to the case where the group orbit through the
relative equilibrium is isotropic}

We begin the proof of Theorem~\ref{main} by reducing it to the case
where the image $\mu =\Phi (m)$ of the relative equilibrium is fixed
by the coadjoint action of $G$, equivalently, when the orbit $G\cdot
m$ is isotropic.
\begin{notation}
If $P$ is a principal $L$-bundle and $N$ an $L$-manifold, we can form
the quotient $P\times_L N:= (P\times N)/L$, the corresponding
associated bundle.  We denote the equivalence class in $P\times_L N  $
of a point $(p,n)\in P\times N$ by $[p, n]$.  If $N$ is a product
$N_1 \times N_2$, the points  of  $P\times_L (N_1 \times N_2)$ will be
written $[p, n_1, n_2]$.
\end{notation}

Suppose a Lie group $L$ acts properly on a manifold $N$, i.e., the map
$L\times N\to N \times N$, $(a, x) \mapsto (a\cdot x, x)$ is proper.
A {\em slice} for this action at $n\in N$ is a submanifold $S$ of $N$,
$n\in S \subset N$,
such that 
\begin{enumerate}
\item $S$ is invariant under the action of the isotropy group $L_n$ of $n$,
\item  $L\cdot S$ is open in $N$,
\item the map $L\times _{L_n}S \to L\cdot S$, $[l, s]\mapsto l\cdot
s$, is an $L$-equivariant diffeomorphism.
\end{enumerate}
Hence if $S$ is a slice, we have $T_n N = T_nS \oplus T_n (L\cdot n)$
and $L_n$ orbits in $S$ parameterize $L$ orbits near $L\cdot n$.

Let $G$ be a Lie group. A point $\alpha \in \fg^*$ is {\bf split} if
its isotropy Lie algebra $\fg_\alpha$ has a $G_\alpha$ invariant
complement $\fn$ in $\fg$, $\fg = \fg_\alpha \oplus \fn$ ($G_\alpha$
equivariant).  Note that if a $G_\alpha$ invariant inner product
exists on $\fg$ (or, equivalently, on $\fg^*$), then $\alpha $ is
split: we can take $\fn$ to be the perpendicular to $\fg_\alpha$ under
the inner product.

Assume that $\alpha \in \fg^*$ is split.
Since the tangent space to the coadjoint orbit $G\cdot \alpha$ is
naturally isomorphic to the annihilator $\fg_\alpha ^\circ$ of
$\fg_\alpha$ in $\fg^*$, one can show that a $G_\alpha$ invariant
neighborhood $B$ of $\alpha$ in the affine subspace $\alpha +
\fn^\circ$ is a slice at $\alpha$ for the coadjoint action of $G$
\cite{GLS}.  Now if $\Phi :M \to \fg^*$ is an equivariant moment map,
then the equivariance of $\Phi$ implies that $\Phi$ is transversal to
$B$.  Hence $R:= \Phi \inv (B)$ is a submanifold of $M$.

\begin{theorem}[{\rm cf.~\cite{GS2}, Theorem~26.7 and \cite{GLS},
Corollary~2.3.6 }] \label{thm:cross-section} Let $(M, \omega )$ be a
Hamiltonian $G$ space with momentum map $F : M\to {\fg}^*$.  Suppose
$\alpha \in {\fg}^*$ is split, and $\fg = \fg_\alpha \oplus \fn$ is a
corresponding $G_\alpha$-invariant splitting.  

Then for a small enough $G_\alpha $ invariant neighborhood $B$ of
$\alpha $ in $\alpha + \fn^\circ$, the preimage $R = F^{-1}
(B)$ is a {\em symplectic} submanifold of $M$. Moreover,  the action of
$G_\alpha $ on $R$ is Hamiltonian with momentum map $F_R : R \to
\fg_\alpha^*$ being the restriction of $F$ to $R$ followed by the
projection  of $\fg^*$ onto $\fg^*_{\alpha}$.
\end{theorem}

\begin{remark}  
The submanifold $R=F ^{-1} (B)$ is called a {\em symplectic
cross-section}.  It has the property that for $m\in F^{-1} (\alpha)$
the $G_\alpha $ orbit is isotropic in $R$.  Also the cross-section $R$
is the smallest symplectic submanifold of $M$ containing the fiber
$F^{-1} (\alpha)$.
\end{remark}

Let $(M, \omega)$ be a symplectic manifold, $G$ a Lie group acting on
$M$ and preserving the symplectic form.  Assume further that there is
an equivariant moment map $\Phi:M \to \fg^*$ corresponding to the
action.  Let $h$ be a $G$ invariant smooth function on $M$ and suppose
$m\in M$ is a relative equilibrium of the Hamiltonian $h$.  Let $\mu =
\Phi (m)$.  Assume that there exists an $H$-invariant inner product on
$\fg$, where $H$ denotes the isotropy group of $\mu$. Then $\mu$ is
split. Let $\fn$ be the $H$ invariant complement of $\fh$ in $\fg$
(orthogonal to $\fh$ with respect to the inner product), $B$ be a $H$
invariant neighborhood of $\mu$ in $\mu + \fn^\circ$ which is a slice
for the action of $G$ and let $R $ be a corresponding cross-section
passing through $m$. Let $\Phi _R : R \to \fh^*$ be a corresponding
moment map.  Since $\mu$ is fixed by $H$, $\pi (\mu) \in \fh^*$ is
fixed by the coadjoint action of $H$ (here $\pi : \fg^* \to \fh^*$ is
the projection).

Let $h_R$ be the restriction of our Hamiltonian $h\in C^\infty (M)^G$
to the cross-section $R$.  Since the flow of the Hamiltonian vector
field $X_h$ of $h$ preserves the fibers of the moment map and the
cross-section is a union of fibers, the vector field $X_h$ is tangent
to the cross-section.  It follows that the Hamiltonian vector field of
$h_R$ on the symplectic manifold $(R, \omega |_R)$ is the restriction
of $X_h$ to $R$.

\begin{lemma} Under the hypothesis of the previous two paragraphs,  if
$m\in R$ is an $H$ stable equilibrium of $h_R$ then $m\in M$ is an $H$
stable equilibrium  of $h$.
\end{lemma}

\begin{proof}
Since $B$ is a slice for the coadjoint action of $G$, $G \cdot B$ is a
neighborhood of $\mu $ in $\fg^*$ diffeomorphic to the associated
bundle $G\times _H B$.  By the equivariance of the moment map, $G\cdot
R$ is an open $G$ invariant subset of the manifold $M$ diffeomorphic
to the associated bundle $G\times _H R$.  The lemma now follows from
the proposition below.
\end{proof}

\begin{proposition}
Let $G$ be a group, $H\subset G$ a compact subgroup. Let $\pi :P=
G\times _H F \to G/H$ be a $G$ homogeneous fiber bundle for some
$H$-manifold $F$.  Suppose $X$ is a $G$ invariant vector field on $P$
which is tangent to the fibers.  Suppose further that $p\in F := \pi
\inv (1H)$ is $H$ stable for the vector field $X$ on $F$ ($1$ denotes the
identity in $G$).  The $p$ is $H$ stable for the vector field $X$ on
the whole space $P$, where $H$ acts on $P$ as a subgroup of $G$.
\end{proposition}

\begin{proof}
Since $X$ is $G$ invariant and is tangent to fibers, its flow $\Phi_t$
is $G$ equivariant and preserves fibers.  Since the action of $G$ on
$P$ is given by $g\cdot [a, f] = [ga, f]$, we have
$$
\Phi_t ([a, f]) = a\cdot \Phi _t ([1, f]).
$$
Since $\Phi_t$ preserves fibers, it induces a flow $\phi_t$ on $F$
by $[1,\phi_t (f)] = \Phi_t ([1, f])$.  Therefore
\begin{equation}
\Phi_t ([a, f]) =a\cdot \Phi_t ([1, f])= =a\cdot [1,\phi_t (f)] =[a,
\phi_t (f)]).
\end{equation}
The subgroup $H$ acts on $G$ by conjugation, $a\cdot g = aga\inv$ for
all $a\in H$, $g\in G$. It also acts on $G/H$ by $a\cdot gH = agH$.
The projection $\varpi: G \to G/H$ is equivariant with respect to these actions.
Since $H$ is compact, there exists a local equivariant section $s$
defined on some $H$ invariant neighborhood $\U$ of $1H$ in $G/H$.  By
definition $s(a\cdot x) = as(x)a\inv$. 
The section defines a diffeomorphism $\psi : \U \times F \to \pi \inv
(\U) \subset P$, $\psi (x, f) = [s(x), f]$.

Since $H$ acts on $F$ and on $\U$, it acts on $\U\times F$: $a\cdot
(x, f) = (a\cdot x, a\cdot f)$.  We have a flow on $\U\times F$:
$\overline{\Phi}_t (x, f) =(x, \phi _t (f))$.  It is not hard to see
that $\psi$ is $H$ equivariant and that it intertwines the flow
$\overline{\Phi}_t$ on $\U\times F$ with the flow $\Phi _t$ on
$\pi\inv (\U)$. Indeed,
\begin{equation*}
\psi ([a\cdot x, a\cdot f]) = [s(a\cdot x),a\cdot  f] = [as(x)a\inv,
a\cdot f] = [as(x), f] = a\cdot [a, f] = a \psi (x, f);
\end{equation*}
and
\begin{equation*}
\psi ([x, \phi _t(f)]) = [s(x),\phi _t(f)] = \Phi_t ([s(x), f]) =
\Phi_t (\psi (x, f)).
\end{equation*}
Consequently it is enough to prove that $(1H, p)$ is $H$ relatively
stable in $\U\times F$ for the flow $\overline{\Phi}_t $.

Let $V$ be an $H$ invariant neighborhood of $(1H, p)$ in $\U\times
F$. Without loss of generality we may assume that $V$ is of the form
$V_1 \times V_2$, where both $V_1 \subset \U$ and $V_2\subset F$ are
$H$ invariant.  Since $p$ is $H$ stable for the flow $\phi_t$, there
is a neighborhood $U_2 \subset V_2$ of $p$ such that $\phi_t
(U_2)\subset V_2$.  Consequently $\overline{\Phi}_t (V_1 \times U_2)
\subset V_1 \times V_2$.
\end{proof}

To complete the reduction to the case where the orbit through the
relative equilibrium is isotropic, it remains to prove one more lemma.

\begin{lemma}
Let $(M, \omega)$ be a symplectic manifold, $G$ a Lie group acting on
$M$ and preserving the symplectic form and $\Phi:M \to \fg^*$ an
equivariant moment map corresponding to the action. Let  $m\in M$ be a
point, $R\subset M$ a symplectic cross-section passing through $m$.
Then 
$$
	\ker d\Phi (m) = \ker d (\Phi |_R) (m).
$$
Hence, if $h\in C^\infty (M)^G$ is an invariant Hamiltonian such that 
$d (h- \l \Phi, \xi\r) (m) = 0$ for some $\xi$ in the isotropy Lie
algebra of $\Phi (m)$ (cf.\ Remark~\ref{rk_isotropy}) , then
$$
  d^2 (h- \l \Phi, \xi\r) (m) |_{\ker d\Phi (m)} = 
d^2 (h|_R- \l  \Phi|_R, \xi\r) (m) |_{ \ker d (\Phi |_R) (m)}
$$
(cf.~Remark~\ref{remark_hessian}).
\end{lemma}

\begin{proof}
The set $G\cdot R$ is a neighborhood of the point $m$ in $M$ and is
diffeomorphic to the associated bundle $G\times _H R$, where $H$ is
the isotropy group of $\mu = \Phi (m)$.  The diffeomorphism $G\times
_H R \to G\cdot R$ is given by $[g, r] \mapsto g\cdot r$.  Similarly,
the set $G\cdot B$ is a neighborhood of the point $\mu$ in $\fg^*$ and
is diffeomorphic to the associated bundle $G\times _H B$, where $B =
\Phi (R)$ (cf.\ Theorem~\ref{thm:cross-section}).  If $s: G/H \supset
U\to G$ is a local section, it simultaneously trivializes the
associated bundles $G\times _H R\to G/H$ and $G\times _H B\to G/H$.
With respect to these trivializations, the moment map $\Phi : M
\supset G\times _H R \supset U\times R \to U\times B \subset G\times
_H B \subset \fg^*$ takes the form
$$
	\Phi (u, r) = (u, \Phi  (r)),
$$
(cf.\ \cite{L_JDG}, p.~814).  Hence for any for any point $r$ in the
cross-section, we have $\ker d\Phi (r) = \ker d (\Phi |_R) (r)$.
\end{proof}

\section{Relative equilibria lying on isotropic orbits}

We now consider the following situation.
Let $(M, \omega)$ be a symplectic manifold, let $G$ be a Lie group
acting properly on $(M,\omega)$ in Hamiltonian fashion and let $\Phi:M \to
\fg^*$ be a corresponding equivariant moment map.  Let $h$ be a $G$
invariant smooth function on $M$.  Suppose that $m\in M$ is a relative
equilibrium of the Hamiltonian $h$.  Suppose further that $\Phi (m)$
is fixed by the coadjoint action of $G$.  Then it is no loss of
generality to assume that $\Phi (m) =0$.  We may also assume that
$h(m) = 0$.

Since the point $0$ is fixed by the action of $G$, the $G$-orbit
through the point $m$ is {\em isotropic} in $(M, \omega)$.  Since the
orbit $G\cdot m$ is isotropic, the tangent space $T _m(G \cdot m)$ is
contained in its symplectic perpendicular $T _m(G \cdot m)^{\omega}$.
Hence the quotient $V = T _m(G \cdot m)^{\omega}/ T_m(G \cdot m)$ is
the {\em symplectic slice \/} at $m$ (cf.\ Remark~\ref{sympl_slice}).
Note that $V$ is a natural symplectic representation of the isotropy
group $K$ of $m$. Note also that since the action of $G$ is proper,
the isotropy group $K$ is compact.

\begin{theorem}[Local normal form for a neighborhood of an isotropic
orbit \cite{GS1}, \cite{marle:model}]
\label{locsymp} 
Let $\Psi : N \to \fh^*$ be a moment map associated to a Hamiltonian
proper action of a  Lie group $H$ on a symplectic manifold
$(N,\sigma)$.  Suppose the orbit $H\cdot x$ is {\em isotropic} in $N$.
Let $H_x$ denote the stabilizer of $x$ in $H$,
let $\fh_x^0$ denote the annihilator of its Lie algebra in $\fh^*$,
and let $H_x \to Sp (V)$ denote the symplectic slice representation.

Given an $H_x$-equivariant embedding, $i:\fh_x^* \to \fh^*$,  
there exists an $H$-invariant symplectic two-form, $\omega_Y$, on the manifold
$Y= H\times_{H_x} \left(\fh_x^\circ \times V \right)$, such that
\begin{enumerate}
\item the form $\omega_Y$ is nondegenerate near the zero section of the
bundle $Y\to H/H_x$, 
\item \label{iso} there exists a neighborhood $U_x$ of the orbit of
$x$ in $N$ and an equivariantly diffeomorphism $\psi$ from a
neighborhood $U_0$ of the zero section in $Y$ to $U_x$ such that
$\psi ^*\omega  = \omega_Y$ , and
\item 
the action of $H$ on $(Y, \omega _Y)$ is Hamiltonian with a moment map
$\Phi_Y : Y \to \fh^*$ given by
$$ 
\Phi _Y ([g, \eta, v])= g\cdot \left(\eta + i( \Phi_{V} (v))\right)
$$ 
where ``$\cdot$'' denotes the coadjoint action, and $\Phi_{V} : V\to \fh_x^*$
is the moment map for the slice representation of $H_x$.
\end{enumerate}
Consequently, the equivariant symplectic embedding $\psi: U_0
\hookrightarrow M$ provided by (\ref{iso}) above intertwines the two
moment maps, up to translation: $\psi^*\Psi = \Phi_Y + \Psi (x)$.
\end{theorem}
\begin{remark}\label{remark_11}
Since the action of $H$ is assumed to be proper, the isotropy group
$H_x$ is compact.  Therefore an $H_x$ equivariant embedding $i:\fh_x^*
\to \fh^*$ always exists: Choose an $H_x$ invariant inner product on
$\fh$.  The orthogonal projection $\fh \to \fh_x$ with respect to the
inner product is $H_x$ equivariant.  Its transpose $i:\fh_x^* \to
\fh^*$ is an $H_x$ equivariant embedding.  Moreover, the perpendicular
$\fh_x^\perp$ of $\fh_x$ (with respect to the inner product) is the
annihilator of $i(\fh_x^*)$ in $\fh^*$.
\end{remark}
\begin{remark}
The form $\omega _Y$ on $Y$ has the further property that orbit
$H\cdot [1,0,0]$ is isotropic and that the symplectic slice at
$[1,0,0]$ is $V$.
\end{remark}
Theorem~\ref{locsymp} is essentially an equivariant version of Weinstein's
isotropic embedding theorem.  The proof in the case that the group $H$ is
compact can be found in  \cite{GS1}  and in \cite{GS2}.  The case of
proper actions is discussed in \cite{LB}.  See also \cite{marle:model}.

Let us now go back to studying the relative equilibrium $m\in M$ lying
on the zero level set of the moment map $\Phi$.  By
Theorem~\ref{locsymp} there exists a $G$-invariant neighborhood $U$ of
$m$ in $M$, a $G$ invariant neighborhood $U_0$ of $[1,0,0]$ in $Y:=
G\times _K (\fk^\circ \times V)$ and an equivariant diffeomorphism
$\psi: U_0 \to U$ such that $\psi ([1,0,0]) = m$, $\psi^* \omega =
\omega _Y$ and $\psi^* \Phi = \Phi_Y$ (since $\Phi _Y ([1,0,0 ]) =
0= \Phi (m) $). Let $h_Y = \psi^* h$.  Then
$$
d\psi (\ker d\Phi _Y ([1,0,0])) = \ker d\Phi (m);  
$$
$$
d (h- \l \Phi, \xi\r) (m) = 0 \quad \text{ if and only if } \quad
d (h_Y- \l \Phi_Y, \xi\r) ([1,0,0]) = 0;
$$
and consequently
$$
\psi^* \left( d^2 (h- \l \Phi, \xi\r) (m)|_{\ker d\Phi (m)} \right) = 
d^2 (h_Y- \l \Phi_Y, \xi\r) ([1,0,0])|_{\ker d\Phi _Y ([1,0,0])} 
$$
(cf.\ Remark~\ref{remark_hessian}).  Moreover if $d^2 (h- \l \Phi,
\xi\r) (m)|_{\ker d\Phi (m)}$ descends to a definite quadratic form on
$V$ then $d^2 (h_Y- \l \Phi_Y, \xi\r) ([1,0,0])|_{\ker d\Phi _Y
([1,0,0])}$ descends to a definite quadratic form on $V$ as well ($V$
is embedded in $Y$ as the subset $\{[1,0,v] \mid v\in V\}$).
Therefore, it is no loss of generality to assume that $(M, \omega) =
(Y, \omega _Y)$ and $m=[1,0,0]$.

Now, if $\xi \in \fk^\perp\subset \fg$, then $\l \Phi, \xi\r ([1,0,
v]) = \l i(\Phi_V (v), \xi \r = 0$ by Remark~\ref{remark_11}.
Therefore $0= d (h- \l \Phi, \xi\r) (m)|_V = d ((h- \l \Phi,
\xi\r)|_V) (m)= d (h |_V) (m)$ and $d^2 (h- \l \Phi, \xi\r) (m)|_V =
d^2 ((h- \l \Phi, \xi\r) |_V) (m) = d^2 (h|_V) (m)$.

Next, observe that the spaces $C^\infty (Y)^G$ and $C^\infty
(\fk^\circ \times V)^K$ are isomorphic: the isomorphism sends a
function $f\in C^\infty (Y)^G$ to the function $\bar{f} \in C^\infty
(\fk^\circ \times V)^K$ defined by $\bar{f}(\eta, v) = f([1,\eta,
v])$.  Therefore the conditions on our Hamiltonian $h\in C^\infty
(Y)^G$ translate into the following two conditions on the
corresponding function $\bar{h}\in C^\infty (\fk^\circ \times V)^K $:
\begin{enumerate}
\item $d_v \bar{h} (0,0) = 0$ and

\item  $d^2_v \bar{h} (0,0)$ is definite.
\end{enumerate}
Here $d_v$ denotes the differential of a function on $\fk^\circ \times
V$ in the $V$ variables and $d^2_v$ has  similar meaning. It is no
loss of generality to assume that the quadratic form  $d^2_v \bar{h}
(0,0)$ on $V$ is positive definite.

Note that the norm squared of the moment map $|\Phi|^2$ is
$G$-invariant and satisfies
$$
|\Phi|^2 ([a, \eta, v]) = |\Phi|^2 ([1, \eta, v]) = |\eta + i(\Phi_V
(v))|^2 = |\eta |^2 +  |i(\Phi_V (v))|^2.
$$
Therefore, if $\gamma (t) = [a(t), \eta (t), v(t)]$ is an integral
curve of the Hamiltonian vector field of $h$ we have
\begin{gather}
|\Phi|^2 (\gamma (0)) = |\Phi|^2 (\gamma (t))\geq |\eta (t) |^2 \\
\text{and} \quad h(\gamma (0))=h(\gamma (t))
\end{gather}
for all time $t$.  Thus the proof of stability reduces to the proof of
the proposition below.

\begin{proposition}
Let $W$ and $V$ be normed finite dimensional vector spaces, $h\in
C^\infty (W\times V)$ a function with $d_v h (0,0) = 0$ and $d^2_v h
(0,0) > 0$, and $\varphi \in C^\infty (W\times V)$ a function with
$\varphi (\lambda, v) \geq |\lambda |^2$.

Then for any $\epsilon >0$ there is $\delta >0$ such that for any
curve $\gamma (t) = (\lambda (t), v(t))$ in $W\times V$ satisfying
$|\gamma (0)|^2 < \delta $, $h(\gamma (t) ) = h (\gamma (0)) $ and
$\varphi (\gamma (t)) = \varphi (\gamma ((0))$ for all $t$, we have
$|\gamma (t)|^2 < \epsilon$ for all $t$.
\end{proposition}

\begin{proof}
We first apply the Morse lemma (see, for example \cite{Duis},
Lemma~1.2.2).  
Since $(0,0)$ is a critical point of $h$ and since the Hessian $d^2 h
(0,0)$ is nondegenerate, there is, by the implicit function theorem, a
function $\sigma (\lambda)$ defined on a neighborhood of $0$ in $W$
such that $d_v h (\lambda, \sigma (\lambda )) = 0$.

Moreover, there exist neighborhoods $U_1$ of $0$ in $W$ and $U_2$ of
$0$ in $V$ and a map $\tau: U_1 \times U_2 \to W\times V$, $\tau
(\lambda, v) = (\lambda, y(\lambda, v))$, $\tau (0,0) = (0,0)$ such
that $\tau (U_1 \times U_2) $ is open $\tau : U_1 \times U_2 \to \tau
(U_1 \times U_2 )$ is a diffeomorphism and
\begin{equation}\label{3.3}
 h(\lam, v) = h( \lam, \sigma (\lam)) + \f{1}{2} d^2_v h (\lam,
\sigma (\lam))(y(\lam, v), y(\lam, v)) .
\end{equation} 

Since $d^2_v h (0,0)$ is positive definite there is $c>0$ such that
$d^2_v h (0,0) (v, v) > c |v|^2$ for all $v\in V$. By continuity,
there is $\epsilon _3 > 0$ such that if $|\lam | < \epsilon_3$ then
$d^2_v h (\lam,\sigma (\lam)) (y, y) > \f{1}{2} c |y|^2$ for all$y\in
V$.

Since $\tau$ is a diffeomorphism onto its image and $\tau (0,0) =
(0,0)$, $\tau \inv$ is continuous near $(0,0)$. Hence for any
$\epsilon >0$ there are $\delta _1, \delta _2 >0$ with
\begin{equation}\label{star}
|\lam|^2 < \delta _1, |y(\lam, v)|^2 < \delta _2 \quad \text{implies }
\quad |\lam|^2 + |y(\lam, v)|^2 < \epsilon.
\end{equation}
Choose $\epsilon_1 >0$ such that $8\epsilon _1/c < \delta
_2$.  Since $\lam \mapsto h(\lam, \sigma (\lam))$ is continuous, there
is $\epsilon _2 >0$ such that 
$$
|\lam|^2 < \epsilon _2 \quad \text{implies }
\quad |h (\lam, \sigma (\lam))|<\epsilon_1.
$$
Let $\epsilon_4 = \min (\epsilon_3, \epsilon_2, \delta _1)$.  Since $h$
and $\varphi$ are continuous, there is $\delta >0$ such that 
$$
|\gamma (0)|< \delta \quad \text{implies} \quad 
\begin{cases} h(\gamma (0)) < \epsilon _2 \quad
\\
\varphi (\gamma (0)) < \epsilon _4 
\end{cases}.
$$
Since $h$ is constant along $\gamma$ and $\varphi (\gamma (0))
=\varphi (\gamma (t))  \geq |\lam (t)|^2$,
$$
|\gamma (0)|< \delta \quad \text{implies}  \quad
\begin{cases} 
h(\lam (t), \sigma (\lam (t))) < \epsilon _2 \quad \\
|\lam (t)|^2 < \epsilon _4 
\end{cases}.
$$
Since 
$$
|\lam (t)|^2 < \epsilon _4  \quad \text{implies} \quad 
\begin{cases} 
d^2_v h (\lam,\sigma (\lam)) (y, y) \geq \f{1}{2} c |y|^2 \\
h (\lam (t), \sigma (t))< \epsilon _1
\end{cases},
$$
we have
\begin{equation*}
\begin{split}
|y(\gamma (t))|^2 &\leq \f{2}{c} d^2_v h (\lam (t), \sigma (\lam
(t)) (y(\gamma (t)),y(\gamma (t)))\\
& \leq \f{4}{c} \left( |h ((\lam (t), v(t))| + |h ((\lam (t), \sigma
(\lam (t)))|\right) \quad \text{by (\ref{3.3})}\\
&  \leq \f{4}{c} (\epsilon _1 +\epsilon _1) < \delta _2  . 
\end{split}
\end{equation*}
Since we arranged $|\lam (t)| < \delta _1$ we have, by (\ref{star}), 
$$
|\gamma (t)|^2 = |\lam(t)|^2 + |v(t)|^2 < \epsilon \quad \text{for all} \, t .
$$
\end{proof}

\end{document}